\documentclass[12pt]{amsart}
\usepackage{amsmath,amsthm,amssymb,euscript}%
\newtheorem{theorem}{Theorem}
\newtheorem{prop}[theorem]{Proposition}
\newtheorem{lemma}[theorem]{Lemma}
\theoremstyle{definition}
\newtheorem{remark}{Remark}
\newtheorem{definition}{Definition}

\newcommand{\PP}{\mathbb P}
\newcommand{\R}{\mathbb R}
\newcommand{\C}{\mathbb C}
\newcommand{\CP}{\C\PP}
\newcommand{\CH}{\C{\mathrm H}}
\newcommand{\Lie}{\mathcal L}
\newcommand{\HH}{\mathbb H}

\newcommand{\CC}{{\EuScript C}}

\newcommand{\J}{{\EuScript J}}
\newcommand{\V}{\mathcal V}

\newcommand{\w}{\omega}

\newcommand{\fhat}{\hat{f}}

\newcommand{\setU}{\EuScript U}
\newcommand{\Uhat}{\widehat{\setU}}

\newcommand{\eg}{\mathbf e}
\newcommand{\eh}{\hat{\mathbf e}}

\newcommand{\bg}{\mathbf g}
\newcommand{\bh}{ h }
\newcommand{\bw}{\mathbf w}
\newcommand{\bn}{\mathbf n}
\newcommand{\bz}{\mathbf z}
\newcommand{\zz}{\mathbf \zeta}
\newcommand{\mean}{{\mathrm m}}
\newcommand{\rv}{\mathrm v}

\newcommand{\ri}{\mathrm i}
\newcommand{\realpart}{\operatorname{Re}}
\newcommand{\impart}{\operatorname{Im}}

\newcommand{\JJ}{{\mathrm J}} 
\newcommand{\di}{\partial}

\newcommand{\intprod}{\mathbin{\raisebox{.4ex}{\hbox{\vrule height .5pt width
5pt depth 0pt %
         \vrule height 3pt width .5pt depth 0pt}}}}

\renewcommand\&{\wedge}

\begin{document}

\title{Hopf Hypersurfaces of Small Hopf Principal Curvature in $\CH^2$}


\author{Thomas A. Ivey}
\address{Dept. of Mathematics, College of Charleston\\
66 George St., Charleston SC 29424-0001}
\email{IveyT@cofc.edu}
\author{Patrick J. Ryan}
\address{Department of Mathematics and Statistics\\
              McMaster University\\
              Hamilton, ON  Canada L8S 4K1}
\email{ryanpj@mcmaster.ca}

\maketitle
\begin{abstract}
Using the methods of moving frames and exterior differential
systems, we show that there exist Hopf hypersurfaces in complex
hyperbolic space $\CH^2$ with any specified value of the Hopf
principal curvature $\alpha$ less than or equal to the
corresponding value for the horosphere.  We give a
construction for all such hypersurfaces in terms of Weierstrass-type data, and also obtain a
classification of pseudo-Einstein hypersurfaces in $\CH^2$.
\end{abstract}

\section{Introduction}
\label{intro} In this paper, we address two classic problems
in the study of hypersurfaces in the
complex space forms $\CP^n$ and $\CH^n$, namely
\begin{itemize}
\item the basic structure theory for Hopf hypersurfaces;
\item the classification of pseudo-Einstein hypersurfaces.
\end{itemize}

\subsection{History of the first problem}  A Hopf hypersurface is one for which the structure vector $W$ is
a principal vector.  (Definitions will be given in \S\ref{definitions}.)
In 1982, Cecil and Ryan \cite{cecilryan} gave a local characterization of Hopf hypersurfaces in $\CP^n$
as tubes over complex submanifolds.  In 1985, Montiel \cite{montiel} established a similar characterization for Hopf
hypersurfaces in $\CH^n$ whose Hopf principal curvature $\alpha$ is greater than $2/r$, the Hopf principal
curvature of the horosphere.

In this paper, we construct the class of Hopf hypersurfaces in $\CH^2$ for which $0 \le \alpha \le 2/r$,
and show that each such hypersurface (for $\alpha <2/r$) can be characterized in terms of Weierstrass-type
data which take the form of a pair of contact curves in $S^3$.  (For a precise statement, see Theorem \ref{mainresult}.)
We expect that a similar approach will show how to construct new examples in higher dimensions.

\subsection{History of the second problem}
The complex space forms $\CP^n$ and $\CH^n$ do not admit Einstein hypersurfaces.  However, as Kon \cite{kon}
discovered, there is a nice class of hypersurfaces satisfying
$$S X = \rho X + \sigma\langle X,W\rangle W,$$
for all tangent vectors $X$, where $S$ is the Ricci tensor, and
$\rho$ and $\sigma$ are constants. He called such hypersurfaces
{\it pseudo-Einstein}.  For $n\ge 3$, the pseudo-Einstein
hypersurfaces were classified by Kon for $\CP^n$ and by Montiel
for $\CH^n$.  Each is an open subset of a homogeneous Hopf
hypersurface (see Theorems 6.1 and 6.2 in \cite{nrsurvey}).  For
$n=2$, the classification problem remained open until recently,
when Kim and Ryan \cite {kimryan2} showed that pseudo-Einstein
hypersurfaces in $\CP^2$ and $\CH^2$ must be Hopf.  However, in
addition to the known homogeneous examples, this classification
includes all Hopf hypersurfaces with $\alpha = 0$.

In this paper, we complete the classification of pseudo-Einstein hypersurfaces in $\CH^2$
by applying our solution to the Hopf hypersurface construction problem for $\alpha=0$.  For more details,
see \S\ref{pseudoEinstein}.

\subsection{Hypersurfaces in Complex Space Forms: Notation and Definitions}\label{definitions}
The complex space forms $\CP^n$ and $\CH^n$ can both be viewed as
quotients of spheres by $S^1$ actions.  In the case of $\CP^n$,
the circle acts on the sphere $S^{2n+1}$ of radius $r$ in
$\C^{n+1}$ by multiplication by a unit modulus complex number. For
$\CH^n$, the action is the same, but the `sphere' is the anti-de Sitter space
$H_1^{2n+1}$ defined by
$\langle \bz, \bz \rangle_\C = -r^2$
using an
Hermitian inner product with signature $(n,1)$:
\begin{equation}\label{hermit}
\langle \bz, \bw \rangle_\C =
-z_0 \overline{w}_0 + \sum_{k=1}^n z_k \overline{w}_k,
\quad \bz,\bw \in \C^{n+1}.
\end{equation}
(In what follows, $\langle\ ,\ \rangle$ without the subscript $\C$ will
indicate the real part of this inner product.  A nonzero vector $\bz$ such
that $\langle \bz, \bz \rangle_\C = 0$ will be called a {\em null vector}.
Also, as in \cite{nrsurvey},
we will use $\HH$ as an abbreviation for $H_1^{2n+1}$.)
In each case, the projection $\pi$ to the quotient is a Riemannian
(respectively, semi-Riemannian) submersion, and
the quotient is a complex space form, i.e., it
 has a positive definite K\"ahler metric of constant holomorphic sectional curvature, given
 by $4/r^2$ or $-4/r^2$, respectively.

A real hypersurface $M$ in $\CP^n$ or $\CH^n$ inherits two
structures from the ambient space.  First, given a unit normal
$\xi$, the {\em structure vector field} $W$ on $M$ is defined so
that
$$\JJ W = \xi,$$
where $\JJ$ is the complex structure.
This gives an orthogonal splitting of each tangent space as
$$ \operatorname{span} \{W\} \oplus W^\perp.$$
Second, we define on $M$ the $(1, 1)$ tensor field $\varphi$ which is the complex
structure $\JJ$ followed by tangential projection, so that
$$\varphi X = \JJ X - \langle X,W\rangle \xi.$$
The shape
operator $A$ is defined by
\begin{equation}
A X = -\nabla_X \xi \label{shape}
\end{equation}
where $\nabla$ is the Levi-Civita connection of the ambient space.
The Gauss equation expresses the curvature tensor of $M$ in terms
of $A$ and $\varphi$ as follows:
$$R(X,Y)=AX\& AY+c(X\& Y+\varphi X\& \varphi Y +2 \langle X,\varphi Y\rangle
\varphi),$$
where $4c = \pm 4/r^2$ is the holomorphic sectional curvature of the complex space form,
and $(X \wedge Y)Z = \langle Y,Z\rangle X - \langle X,Z\rangle Y$.

The Ricci tensor is given by
\begin{equation}\label{Ricci}
SX=(2n+1)cX-3c\langle X,W\rangle W+ \mean AX-A^2 X,
\end{equation}
 where $\mean = {\text{trace}}\ A$.

Let $\alpha = \langle AW, W \rangle$.  If $W$ is a principal vector everywhere
(i.e. $A W = \alpha W$), we say that $M$ is a {\it Hopf hypersurface}.
Following J.K. Martins \cite{martins}, we use the term {\em Hopf
principal curvature} to mean the principal curvature that
corresponds to the principal vector $W$ for a Hopf hypersurface.

It is a nontrivial fact that the Hopf principal curvature $\alpha$ is
constant for a (connected) Hopf hypersurface.
This was proved by Y. Maeda \cite{YMaeda} for $\CP^n$
and by Ki and Suh \cite{KiSuh} for $\CH^n$.
Further details on the material of this section may be found in \S2 of \cite{nrsurvey}.

\subsection{Our approach}
Many classification results in the study of hypersurfaces
in complex space forms rely on the ingenious
use of the identities that arise between the shape operator, the
curvature tensor, and $\varphi$.  Unfortunately, some of these
arguments work only when the complex dimension $n$ of the space
form is at least three.  Motivated by some unresolved questions in
the case $n=2$, we began a program of using moving frames and
exterior differential systems to study special hypersurfaces in
$\CP^2$ and $\CH^2$.
This approach has the advantage of being
systematic.  One does not have to choose which tensor to
differentiate in which direction; instead, Cartan's test for
involutivity, alternating with prolongation, is guaranteed to tell
you how many solutions there are.  Nonetheless, some effort
is required to set up the machinery of the frame bundle and its
canonical framing, in terms of which our exterior differential
systems are defined.  Although the results of this paper do not
make use of Cartan's test, we used it to show existence of Hopf
hypersurfaces with $0 \le |\alpha| \le 2/r$ after which we were
able to make an explicit construction as described in Theorem \ref{mainresult}.

In this paper, all manifolds are assumed connected and all
manifolds and maps are assumed smooth $(C^\infty)$ unless stated otherwise.  For basic
reference material on hypersurfaces, see \cite{nrsurvey}.  More on
exterior differential systems may be found in the monograph
\cite{BCG3} or the textbook \cite{cfb}.

\section{Setup and Statement of Results}
In this section we state our main result (Theorem \ref{mainresult} below).  In order
to make sense of the mappings used in that theorem, we first need to define moving frames
for hypersurfaces in $\CH^2$.  In fact, we will work with lifted frames, which are
sections of the frame bundle of $\HH = H_1^5$ that are orthonormal with respect to
the semi-Riemannian metric given by the real part of \eqref{hermit}.
Because this may be less familiar to some readers than the Riemannian case,
and also in order to fix notation, we review the basic tools in some detail.

\subsection{Moving Frames for Hypersurfaces in $\CH^2$}\label{movingframes}

Let $G$ be the Lie group $U(2,1)$ of matrices whose columns form a `unitary' basis
for $\C^3$ with respect to the Hermitian inner product \eqref{hermit} with signature $(2,1)$.
We take the convention that the {\em first} column of a matrix $u\in G$ has inner product $-1$
with itself, while the second and third columns have inner product $+1$ with themselves.

Let $\zz:G\to \HH$ be given by multiplying the first column of $u \in G$ by $r$,
and let $\rho = \pi \circ \zz: G \to \CH^2$.
Let $\eg_2$ and $\eg_3$ denote the second and third columns
of $u$.
At the point $\bz=\zz(u) \in \HH$, vectors $\eg_2(u),\eg_3(u)$ span the complex plane that
is orthogonal to $\bz$; thus, these vectors push forward under $\pi_*$
to span the tangent space of $\CH^2$ at $\pi(\bz)$.
For later use, we will need to consider frames that span these spaces
as real vector spaces.
To that end, define additional $\C^3$-valued functions $\eg_0,\eg_1,\eg_4$ on $G$ such
that
\begin{equation}\label{ehatrels}
\dfrac{\ri}{r} \zz = \eg_0,\qquad \ri \eg_1 = \eg_2, \qquad \ri\eg_3 =\eg_4.
\end{equation}
The vector $\eg_0(u)$ is tangent to the fiber of $\pi$ at $\bz$, normalized so that
$\langle \eg_0, \eg_0 \rangle_\C=-1$, while the remaining
unit vectors $\eg_1(u),\ldots, \eg_4(u)$ push forward under $\pi$ to give a basis (over $\R$) for the
tangent space to $\CH^2$ at $\pi(\bz)$.  (These are referred to as {\em horizontal} vectors, as
they are orthogonal to the fiber of $\pi$.)

\begin{definition}
Let $M^3 \subset \CH^2$ be an oriented embedded hypersurface.
We say that a map $f:M \to G$ is an {\em adapted lift} of $M$ if $\rho \circ f = \operatorname{id}_M$,
$\pi_* (\eg_4 \circ f)$ is normal to $M$ and agrees with the orientation.
(Note that this means that $\pi_*(\eg_3 \circ f)$ is the structure vector $W$ of $M$ and that
$\pi_* (\eg_1 \circ f)$ and $\pi_*(\eg_2 \circ f)$ are also tangent to $M$.)
\end{definition}
The choice of adapted lift is not unique; at a given
point of $M$, we may change $\bz$ by the $S^1$ action on $\HH$, and we may also multiply $\eg_2$
by a unit modulus complex number.

\subsection{Gauss Maps and Characteristics}\label{gausschar}
Now we will specialize to Hopf hypersurfaces $M\subset \CH^2$ for a fixed Hopf principal
curvature $\alpha \in (-2/r, 2/r)$.  (We will discuss the case where $|\alpha|=2/r$ in \S\ref{sideproof}.)
Write $\alpha = \frac{2}{r}\sin\phi$ where $\phi \in (-\pi/2,\pi/2)$.
We define two additional vector-valued functions on $G$,
\begin{equation}\label{gdef}
\bg^\pm = \eg_0 -(\sin\phi\, \eg_3 \pm \cos\phi\, \eg_4).
\end{equation}
and let $\bg^\pm_\C = \pi \circ \bg^\pm$ denote the line in $\C^3$ spanned by $\bg^\pm$.

For an adapted lift $f:M \to G$, the compositions $ \bg^\pm\circ f$
may be regarded as modified versions of the
Gauss map:  first, the hypersurface normal is rotated within a complex plane,
through an angle determined by $\alpha$, and then lifted
to a null vector tangent to the anti-de Sitter space $\HH$.

Let $\V$ denote the set of nonzero null vectors in $\C^3$.  Then the $\bg^\pm$ take values in $\V$, and
the $\bg^\pm_\C$ are maps from $G$ into the projectivized null cone
$\pi(\V) \subset \CP^2$.
This space is a smooth manifold of real dimension three, and we will identify it with $S^3$ as follows:
Given a point $\bz \in \V$, satisfying
\begin{equation}\label{znull}
|z_1|^2 + |z_2|^2 = |z_0|^2,
\end{equation}
we note that $z_0 \ne 0$, and thus $\pi(\V)$ lies entirely in the domain of one of the
standard coordinate charts on $\CP^2$.  Letting $w_1=z_1/z_0$ and $w_2=z_2/z_0$, we see immediately
that points in $\pi(\V)$ satisfy $|w_1|^2 + |w_2|^2=1$, the equation of the unit sphere in $\C^2$.

The relationship between the geometry of the Hopf hypersurface $M$
and the null vectors $\bg^\pm$ is as follows.
Choose any point $p \in M$ and an orthonormal principal basis $(e_1, e_2, e_3)$ for $T_p M$
satisfying
$e_2 = \varphi e_1$ and $e_3 = W$, so that
$A e_1 = \lambda e_1$, $Ae_2 = \nu e_2$, and $A e_3 = \alpha e_3$.
A fundamental identity for Hopf hypersurfaces
(following from the Codazzi equation -- see \cite{nrsurvey}, pp. 245--246) gives
\begin{equation} \label{Hopfpc}
\lambda \nu = \frac {\lambda + \nu}{2}\, \alpha + c
\end{equation}
(recall that $c=-1/r^2$).
This shows that $\lambda$ and $\nu$ are distinct, since
\begin{equation}\label{splitsig}
\left(\lambda-\dfrac{\alpha}2\right)\left(\nu- \dfrac{\alpha}2\right)
= c +\dfrac{\alpha^2}4 = -\dfrac{\cos^2\phi}{r^2}.
\end{equation}
Hence $(e_1, e_2, e_3)$ extends to a smooth local principal orthonormal
frame near $p$ with corresponding smooth principal
curvatures $(\lambda, \nu, \alpha)$.
Moreover, \eqref{splitsig} implies that the quadratic form
corresponding to the restriction of $A-\frac{\alpha}2 I $ to $W^\perp$
has signature $+-$.
We define two vectors $\rv^\pm$ in $W^\perp$ that are null for this quadratic form, i.e.
\begin {equation}
 \langle (A - \frac{\alpha}{2})\rv^{+} , \rv^{+} \rangle = 0, \quad
\langle (A - \frac{\alpha}{2})\rv^{-} , \rv^{-} \rangle = 0,
\end{equation}
given by
\begin{equation}\label{asymptotics}
\rv^\pm = \pm \cos \phi\ e_1 + (r \lambda -\sin \phi)\varphi e_1.
\end{equation}
Although the values of $\rv^+$ and $\rv^-$ depend on the choice of
which eigenvector is assigned to $e_1$, it is easy to check
using \eqref{Hopfpc} that a different choice merely multiplies
each of $\rv^\pm$ by a nonzero factor. Thus, the respective distributions spanned by $\rv^+$
and $\rv^-$ are well-defined and $M$ is foliated in two ways by curves
tangent to these distributions.

We define the {\em characteristic distributions} $\chi^+$ and $\chi^-$ on $M$
as the direct sum of the span of $\rv^+$ (respectively, $\rv^-$) and the span of the
structure vector $W$.
It is somewhat surprising that these two-dimensional distributions on $M$ are
integrable; in fact, we have
\begin{prop}\label{intprop}  Let $M \subset \CH^2$ be a Hopf hypersurface with
Hopf principal curvature $\alpha \in (-2/r,2/r)$.  Then $M$ is foliated by
two-dimensional leaves of each of the distributions $\chi^+, \chi^-$.  Moreover,
if $f$ is an adapted lift of $M$, then the maps $\bg^+_\C \circ f$ and $\bg^-_\C \circ f$ are constant
along the leaves of $\chi^+,\chi^-$, respectively.
\end{prop}
This proposition will be proved in \S\ref{charproof}.

\subsection{Main Result}
\begin{theorem}\label{mainresult}\
\begin{enumerate}
\item Let $M$ be an oriented Hopf hypersurface in $\CH^2$
with Hopf principal curvature $\alpha \in (-2/r, 2/r)$.  Then there are well-defined
maps $\sigma^\pm:M \to S^3$, such that $\sigma^\pm = \bg^\pm_\C \circ f$ for any adapted
lift $f$ of $M$, and such that the images of $\sigma^\pm(M)$ are contact curves\footnote{The contact
structure here is the standard one in $S^3$, whose contact planes are orthogonal to the fibers of
the Hopf fibration $S^3 \to \CP^1$.}
 in $S^3$.

\item Conversely, given any two embedded contact curves $\CC_1, \CC_2$ in $S^3$, let $P$ be
the intersection of their inverse images under $\bg^\pm_\C$ in $G$.
Then $\rho(P)$ is a Hopf hypersurface with Hopf principal curvature $\alpha$.
\end{enumerate}
\end{theorem}

We will prove this theorem in \S \ref{mainproof}.  For the moment, we note how these
results carry over to the borderline case, when $\alpha = \pm 2/r$.
In this case, $\cos\phi=0$, so $\bg^+=\bg^-$.  Labelling this as $\bg$,
we note that $\bg \circ f$ is well-defined (in fact, it is $\eg_0-\eg_3$, which is a just
a null lift of $-W$, i.e., one of two null vectors tangent to $\HH$ that project to $-W$).
The image of the map $\bg_\C$ is still a contact curve in $S^3$,
which will be regular at points where $\lambda \ne \nu$.
However, given the image contact curve,
more data is necessary to reconstruct the Hopf hypersurface $M$;
more details will be given in \S\ref{sideproof}.

As for the characteristic distributions, in the borderline case \eqref{splitsig} becomes
\begin{equation}\label{borderline}
\left(\lambda-\dfrac{\alpha}2\right)\left(\nu- \dfrac{\alpha}2\right)= 0.
\end{equation}
If $\lambda \ne \nu$, we may assume, without loss of generality,
that $\lambda \ne \dfrac{\alpha}2$ so that $\rv^+ = \rv^-$ and
the {\it unique} characteristic distribution is spanned by $\varphi e_1$, the principal direction
corresponding to $\nu = 1/r$.
On the other hand, if $\lambda = \nu$ globally, then $M$ is just an open subset of a horosphere with $\lambda = 1/r$.  Finally, there is the possibility that $\lambda = \nu$ at some points and $\lambda \ne \nu$ at others.

\begin{remark}
Tubes over real hyperbolic planes occur in the list of standard
examples of Hopf hypersurfaces
with constant principal curvatures.  As seen in  Theorems 3.4 and 3.12 of \cite{nrsurvey}, they
have Hopf principal curvature $0 < \alpha < 2/r$.  According to
Theorem \ref{mainresult}, these must be constructible by our
methods.  It would be interesting to identify which contact curves produce these special
examples.
\end{remark}

\subsection{Classification of Pseudo-Einstein Hypersurfaces}\label{pseudoEinstein}
As defined by Kon \cite{kon},
a hypersurface $M^{2n-1}$ in a complex space form is {\em pseudo-Einstein} if there
are constants $\rho$ and $\sigma$ such that the Ricci tensor $S$
of $M$ satisfies
$$S X = \rho X + \sigma\langle X,W\rangle W,$$
for all tangent vectors $X$, where $W$ is the structure vector.  When $n\ge 3$, it
is known that every such hypersurface is an open subset of a homogeneous
Hopf hypersurface (see Theorems 6.1 and 6.2 in \cite{nrsurvey}).  When $n=2$,
we have

\begin{theorem}[Kim-Ryan \cite{kimryan2}]
Every pseudo-Einstein hypersurface $M$ in $\CH^2$ or $\CP^2$ is a Hopf hypersurface.
In $\CP^2$, $M$ is either an open subset of a geodesic sphere, or has
$\alpha=0$ and is generically an open subset of a tube of radius $r \pi/4$ over
a holomorphic curve.  In $\CH^2$, either $\alpha = 0$ or $M$ is an open subset of a horosphere or a tube
over a totally geodesic $\CH^1$.
\end{theorem}

Moreover, it is easy to check that every Hopf hypersurface in
$\CH^2$ with $\alpha=0$ is pseudo-Einstein.  In fact, equation (2.2b)
and Proposition 2.21 of
\cite{kimryan2} shows that for such a hypersurface,
$SW = 2c W$ and $SX = 6c X$ for all $X \in W^{\perp}$.
One can also see this from \eqref{Hopfpc} and \eqref{Ricci} in the current paper.
First note that $\lambda \nu = c$ from \eqref{Hopfpc} so that in \eqref{Ricci} we have
$\mean A X -A^2 X = cX$ for $X \in W^{\perp}$ and thus $SX = 6c X$.  The fact that $SW = 2c W$ when
$\alpha = 0$ is immediate.

Thus, the construction given earlier
in this section, in the special case where $\alpha=0$, completes
the classification of pseudo-Einstein hypersurfaces in complex
space forms.  We have

\begin{theorem}
A pseudo-Einstein hypersurface in $\CH^2$ must be a Hopf hypersurface.  Denoting the Hopf principal curvature by
$\alpha$, which we may arrange to be nonnegative, we have the following possibilities:
\begin{enumerate}
\item $\alpha > 2/r$.  An open subset of a geodesic sphere or tube over a totally geodesic $\CH^1$,
as described in Theorems 3.7 and 3.8
of \cite{nrsurvey}.
\item $\alpha = 2/r$.  An open subset of a horosphere, as described in Theorem 3.4 of \cite{nrsurvey}.
\item $\alpha = 0$.  A hypersurface constructed from two arbitrary contact curves as in
Theorem \ref{mainresult} (using $\alpha=0$).
\end{enumerate}
\end{theorem}

\section{Proofs}
In this section we will give proofs of Proposition \ref{intprop} and Theorem \ref{mainresult}.
Because it will be convenient to compute using differential forms, we begin by deriving
the structure equations for a basis for the left-invariant 1-forms on $G=U(2,1)$.

\subsection{Structure Equations}\label{streqsec}

Recall from \S\ref{movingframes}
that we define vector-valued functions $\zz, \eg_2, \eg_3$ on $G$ by decomposing a matrix $u\in G$ into columns:
\begin{equation}\label{u21matrix}
u = [ \dfrac1r \zz(u), \eg_2(u), \eg_3(u)].
\end{equation}
Let $\gamma^i_j$, for $1\le i,j\le 3$, denote the components of
the Maurer-Cartan form $u^{-1} du$ on $G$.  Note that these are
complex-valued 1-forms. Then
\begin{equation} \label{u21streqs}
\begin{aligned}
d\zz &=  \gamma^1_1 \zz + r( \gamma^2_1\eg_2 + \gamma^3_1\eg_3 ),\\
d\eg_2 &= \tfrac1r \gamma^1_2\zz  + \gamma^2_2\eg_2  + \gamma^3_2\eg_3 ,\\
d\eg_3 &= \tfrac1r \gamma^1_3\zz  + \gamma^2_3\eg_2  + \gamma^3_3\eg_3 .
\end{aligned}
\end{equation}
Differentiating the relations $\langle \zz, \zz \rangle_\C=-1/r^2$, $\langle \zz, \eg_2 \rangle_\C=0$,
$\langle \eg_2, \eg_2 \rangle_\C=1$, etc., and using \eqref{u21streqs}, shows that
$$\overline{\gamma^j_i} = \left\{\begin{aligned}
\gamma^i_j &\text{ if exactly one of $i,j$ equals 1}\\
-\gamma^i_j &\text{ otherwise.}\end{aligned}
\right.
$$
As well, the $\gamma^i_j$ satisfy the usual Maurer-Cartan equations:
\begin{equation}\label{gammaMC}
d\gamma^i_j = -\gamma^i_k \& \gamma^k_j.
\end{equation}
(We use the usual summation convention from now on.)

We will identify $G$ with a sub-bundle of the orthonormal frame bundle of $\HH$, by associating to
$u \in G$ the orthonormal frame $(\eg_0(u), \ldots, \eg_4(u))$ at base point $\bz$.
(Recall that $\eg_0,\eg_1,\eg_4$ are related to $\zz,\eg_2,\eg_3$ by \eqref{ehatrels}.)
To analyze the geometry of real hypersurfaces in $\CH^2$, we will need to use
real-valued differential forms on the frame bundle.  To that end,
we first define real-valued 1-forms $\eta^0, \ldots, \eta^4$  on $G$ such that
\begin{equation}\label{dzed}
d\zz = \eta^0 \eg_0 + \eta^1 \eg_1 +\eta^2 \eg_2 + \eta^3 \eg_3 + \eta^4 \eg_4.
\end{equation}
Comparison with \eqref{u21streqs} shows that the $\eta^i$ are
linearly independent left-invariant 1-forms on $G$, and that
\begin{equation}\label{diagamma}
\gamma^1_1 = \tfrac{\ri}r \eta^0, \qquad
\gamma^2_1 = \tfrac1r(\eta^2 - \ri \eta^1), \qquad \gamma^3_1 = \tfrac1r(\eta^3+\ri \eta^4).
\end{equation}
Moreover, \eqref{dzed} shows that the $\eta^i$ are the canonical forms on $G$,
once we identify $G$ with a frame bundle.

The connection forms $\w^\alpha_\beta$ on $G$, where $0 \le \alpha,\beta \le 4$,
are uniquely characterized by satisfying the structure
equations
\begin{equation}\label{characterize}
d\eta^\alpha = -\w^\alpha_\beta \& \eta^\beta
\end{equation}
and taking value in the Lie algebra $\mathfrak{o}(1,4)$, i.e., $\w^j_0 = \w^0_j$ and
$\w^j_i = -\w^i_j$.  (We now take the index ranges $1 \le j,k \le 4$.)
Thus, to calculate the connection forms we need to
express the derivatives of the canonical forms in terms a basis of 1-forms on $G$.
These 1-forms are in turn defined by differentiating the frame vectors, regarded
as vector-valued functions on $G$.

Multiplying \eqref{dzed} by $\ri/r$, we have
\begin{equation}\label{dehatzero}
d\eg_0 = -\frac1{r^2} \eta^0 \zz  + \frac1r  \J^j_k \eta^k\eg_j,
\qquad\text{where}\ \J = \left[\begin{array}{c|c}
\begin{smallmatrix} 0 & -1 \\ 1 & 0 \end{smallmatrix} & 0 \\[2pt] \hline
0 & \begin{smallmatrix} 0 & -1 \\ 1 & 0 \end{smallmatrix}
\end{array}\right].
\end{equation}
The matrix $\J$ represents multiplication by $\ri$ on the span of $\{\eg_1, \ldots, \eg_4\}$.
Differentiating $\langle \zz, \eg_j \rangle$ and $\langle \eg_0, \eg_j \rangle$ shows that
\begin{equation}\label{dhatalpha}
d\eg_j = \eta^k_j\eg_k  + \frac1{r^2} \eta^j\zz
+ \frac1r\J^j_k \eta^k \eg_0,
\end{equation}
for some (real-valued) 1-forms $\eta^j_k$.   By differentiating
$\langle \eg_j, \eg_k \rangle$, we see that the $\eta^j_k$ are
skew-symmetric.  In fact, comparing with \eqref{u21streqs} shows
that
\begin{equation}\label{offdiagamma}
\eta^2_1 = \impart \gamma^2_2, \quad
\eta^3_1 = \eta^4_2 = \impart \gamma^3_2,\quad
\eta^3_2 = -\eta^4_1 = \realpart \gamma^3_2,\quad
\eta^4_3 = \impart \gamma^3_3.
\end{equation}
The 1-forms $\eta^0, \ldots, \eta^4$, $\eta^2_1$, $\eta^4_1$, $\eta^4_2$, $\eta^4_3$
are linearly independent, and form a basis for the left-invariant 1-forms on $G$.

Using \eqref{diagamma}, \eqref{offdiagamma} and the Maurer-Cartan equations \eqref{gammaMC},
we compute that
\begin{align}
d\eta^0 &= -\tfrac1r \J^k_j \eta^j \& \eta^k,\\
d\eta^j &= -(\eta^j_k-\tfrac1r \J^j_k \eta^0 ) \& \eta^k. \label{detaj}
\end{align}
Comparing these equations with \eqref{characterize} shows that
the connection forms are given by $\w^0_k = \w^k_0 = \tfrac1r \J^k_j \eta^j$ and
$$\w^j_k = \eta^j_k -\tfrac1r \J^j_k \eta^0.$$

Next, let $F$ be the orthonormal frame bundle of $\CH^2$.  We define a map $\Pi :G \to F$
which takes $u\in G$ to an orthonormal frame at $\pi(\zz(u))$ given
by $(\pi_* \eg_1, \ldots, \pi_* \eg_4)$.
The relationship between the connection forms on $F$ and on $G$ is given by
the following

\begin{lemma}\label{connectionvlemma}
Let $f_0=(e_1, e_2, e_3, e_4)$ be an orthonormal frame on $U \subset \CH^2$
such that $\JJ e_1 = e_2$ and $\JJ e_3 = e_4$,  and let $f:U \to G$ be any map
such that $\Pi \circ f = f_0$.  Then for any tangent vector $\rv$ on $U$,
\begin{align}
\rv &= (\rv \intprod f^* \eta^k) e_k, \label{omegacanon}\\
\nabla_\rv e_j &= (\rv \intprod f^* \w^k_j ) e_k, \label{omegaconn}
\end{align}
where $\nabla$ denotes the Levi-Civita connection on $\CH^2$.
\end{lemma}
\begin{proof} Let $\bh = \zz \circ f$, and
let $\eh_1, \ldots, \eh_4$ be the horizontal vector fields  on
$\bh(U)\subset \HH$ such that $(\pi_*)_{\bh(p)} \eh_j(\bh(p)) = e_j(p)$
for all $p \in U$. Note
that $\eh_j \circ \bh =\eg_j \circ f$. We may extend these vector
fields to all of $\Uhat = \pi^{-1}(U)$ using the $S^1$-action, and
still have $\pi_* \eh_j = e_j$.  These vector fields, together
with $\eh_0(\bz) = \tfrac{\ri}{r}\bz$, define a section $\fhat:
\Uhat \to G$ such that $f = \fhat \circ \bh$.

\medskip
Suppose $\rv \in T_p \CH^2$, and let $\hat\rv \in T_{\bh(p)} \HH$
be the horizontal vector that projects to $\rv$.
By the defining property of the canonical forms on $G$,
$$\hat\rv = (\hat\rv \intprod \fhat^* \eta^k) \eh_k.$$
Then, because $\hat\rv$ and $\bh_*\rv$
differ by a multiple of $\eh_0$,
$$\hat\rv =((\bh_*\rv) \intprod \fhat^*\eta^k)\eh_k
=(\rv \intprod f^* \eta^k) \eh_k.$$
Pushing forward by $\pi_*$ gives \eqref{omegacanon}.

Let $\widehat\nabla$ be the Levi-Civita connection for the semi-Riemannian metric on $\HH$.
Because $\pi:\HH \to \CH^2$ is a Riemannian submersion,
then $\nabla_{\rv} e_k$ equals the pushforward, via $\pi_*$,
of the horizontal part of $\widehat \nabla_{\hat \rv} \eh_k$.
(For more details, see \S1 of \cite{nrsurvey}.)
Then, by the defining property of the connection forms on $G$,
$$\widehat \nabla_{\hat\rv} \eh_k = (\hat\rv \intprod \fhat^*\w^j_k) \eh_j + (\hat\rv \intprod \fhat^*\w^0_k) \eh_0.$$
Pushing forward by $\pi_*$ gives
$$\nabla_{\rv} e_k = (\hat\rv \intprod \fhat^*\w^j_k) \pi_*\eh_j  = (\rv \intprod f^* \w^j_k) e_j.$$
Note that $\eh_0\intprod\fhat^*\w^j_k = \langle \eh_j, \widehat\nabla_{\eh_0} \eh_k \rangle = 0$
because the $\eh_k$ are translated by parallelism along the fibres of $\pi$.
\end{proof}

Lemma \ref{connectionvlemma} implies that the pullbacks, under
$\Pi:G\to F$, of the canonical and connection forms on $F$ are the
$\eta^k$ and $\w^j_k$ on $G$. Any computation made on $G$ using
these pullback forms will be valid on $F$.

\subsection{Characteristics}\label{charproof}

\begin{lemma}\label{IIlemma}  Let $M \subset \CH^2$ be a hypersurface, and let $f: M \to G$ be
an adapted lift of $M$.  Let $e_i = \pi_*( \eg_i \circ f)$ for $1 \le i \le 4$.
Then $f^*\eta^4=0$ and
\begin{equation}\label{fajk}
f^* \w^4_j = \sum_{k=1}^3 \langle A e_j, e_k\rangle\, f^* \eta^k,
\end{equation}
for $1\le j \le 3$.
\end{lemma}
\begin{proof}Near any point $p\in M$, we can extend the domain of $f$ to
an open set $U$ around $p$.  (Of course, the condition that the vector $\eg_4 \circ f$
projects to be normal to $M$ does not apply at points not on $M$.)
Then, applying \eqref{omegacanon} to vectors $\rv$ tangent to $M$ shows that
$f^*\eta^4=0$.

Furthermore, for $1 \le j, k \le 3$, we have $(f^*\eta^k)(e_j) = \delta^k_j$ on $M$.
Thus, using \eqref{omegaconn}, we have for $1 \le \ell \le 3$,
\begin{align*}
\sum_{k=1}^3 \langle A e_j, e_k\rangle\, f^* \eta^k(e_{\ell}) &=
\sum_{k=1}^3 \langle A e_j, e_k\rangle\ \delta^k_{\ell}= \langle A e_{\ell}, e_j\rangle \\
&=
 -\langle \nabla_{e_{\ell}} e_4, e_j\rangle = - f^* \w^j_4(e_{\ell}) = f^* \w^4_j(e_{\ell})
 \end{align*}
which establishes the desired result.
\end{proof}

\begin{proof}[Proof of Proposition \ref{intprop}]
We will verify the statements in the proposition using
an adapted lift $f$ of the Hopf hypersurface $M$.  (Such adapted
lifts exist locally near any point of $M$.)  Let $e_1, \ldots, e_4$
be as in the proof of Lemma \ref{IIlemma}.  By the Hopf condition,
$$\langle A e_1,\,e_3\rangle = \langle A e_2,\, e_3\rangle =0, \qquad \langle A e_3,\, e_3\rangle =\alpha.$$
Thus, $f^*(\w^4_3 - \alpha  \eta^3)=0$.
(For the rest of this section, differential forms will be understood
to be pulled back by $f$.)

Differentiating $\eta^4=0$ and $\w^4_3 -\alpha \eta^3 =0$ (and substituting for
$\eta^4$ and $\w^4_3$ using these equations) gives
\begin{gather}
\w^4_1 \& \eta^1 + \w^4_2 \& \eta^2 = 0, \label{first2form} \\
2 \w^4_1 \& \w^4_2 + \alpha (\w^4_2 \& \eta^1 - \w^4_1 \& \eta^2) + \frac{2}{r^2} \eta^1 \& \eta^2
=0. \label{second2form}
\end{gather}

Linearly combining \eqref{first2form} and \eqref{second2form} gives two equivalent equations,
\begin{equation}\label{recombined}
\begin{aligned}
(\w^4_1 -\frac1r(\sin\phi\, \eta^1 + \cos\phi\, \eta^2))
\& (\w^4_2 -\frac1r(\sin\phi\, \eta^2 - \cos\phi\, \eta^1))&=0,\\
(\w^4_1 -\frac1r(\sin\phi\, \eta^1 - \cos\phi\, \eta^2))
\& (\w^4_2 -\frac1r(\sin\phi\, \eta^2 + \cos\phi\, \eta^1))&=0.
\end{aligned}
\end{equation}
(Recall that $\alpha = (2/r)\sin\phi$.)
For the sake of convenience, we define
\begin{align*}
\kappa^\pm_1 &= \w^4_1 -\frac1r(\sin\phi\, \eta^1 \pm \cos\phi\, \eta^2),\\
\kappa^\pm_2 &= \w^4_2 -\frac1r(\sin\phi\, \eta^2 \mp \cos\phi\, \eta^1).
\end{align*}
Then the equations \eqref{recombined} imply that $\kappa^+_1$ and $\kappa^+_2$
are linearly dependent, as are $\kappa^-_1$ and $\kappa^-_2$.  (Note, however,
that at least one form in each pair is nonzero at every point.)

Recall the characteristic distributions on $M$ defined in \S\ref{gausschar}.
There is no loss of generality here in assuming that our frame is principal on $M$,
i.e $Ae_1 = \lambda e_1$ and $Ae_2 = \nu e_2$ for distinct principal curvature
functions $\lambda, \nu$.
Using $\w^4_1= \lambda \eta^1$ and $\w^4_2 = \nu\eta^2$,
it is easy to check that $\kappa^+_1, \kappa^+_2$ vanish on
 $\chi^+$, and $\kappa^-_1, \kappa^-_2$ vanish on $\chi^-$.
To show integrability, suppose that $\kappa^+_1 \ne 0$ near
a point $p \in M$. Choose 1-forms $\tau_1, \tau_2$ so that $\kappa^+_1 \& \tau_1 \& \tau_2 \ne 0$.  Then,
$d\kappa^+_1 = \kappa^+_1 \& \theta + k\, \tau_1 \& \tau_2$ for some 1-form $\theta$ and some scalar $k$.
Once we check that $\kappa^+_1 \& d\kappa^+_1=0$, we get $k=0$, which establishes the integrability
condition for $\chi^+$ near $p$.
Using $\eta^4=0$ and $\w^4_3 = \alpha \eta^3$, we calculate
that
\begin{align}
d\kappa^+_1 &=\tfrac1r(\cos\phi\, \kappa^+_1 +\sin\phi\, \kappa^+_2) \& \eta^3 + \eta^2_1 \& \kappa^+_2.
\end{align}
Since $\kappa^+_1$ and $\kappa^+_2$ are linearly dependent, we must have $d\kappa^+_1 \& \kappa^+_1 = 0$.

A similar argument works for $\kappa^+_2 \ne 0$ and for the distribution $\chi^-$. This proves the first assertion
of the proposition.

Next, we show that the maps $\bg^\pm_\C$ are constant along the corresponding characteristic surfaces in $M$.
For simplicity, we check this assertion for $\bg^+_\C$, the argument for $\bg^-_\C$ being similar.  First, before
restricting to $M$, we use the equations \eqref{dehatzero} and \eqref{dhatalpha} to calculate
\begin{multline}\label{dgplus}
d\bg^+ = \frac1r (\sin\phi\, \eta^4 - \cos \phi\, \eta^3)\bg^+
+\frac1r  (\eta^0 + \sin \phi\, \eta^3 + \cos\phi\, \eta^4)\ri\bg^+  \\
-\kappa^+_1(\sin\phi\, \eg_2-\cos\phi\, \eg_1)
+ \kappa^+_2 (\sin\phi\, \eg_1 + \cos\phi\,\eg_2)
-\kappa^+_3 (\sin\phi\, \eg_4 - \cos\phi\,\eg_3),
\end{multline}
where, for future reference, we define
$$\kappa^\pm_3 = \w^4_3 - \tfrac{2}{r}(\sin\phi\, \eta^3 \pm \cos\phi\ \eta^4)
=\w^4_3 -\alpha \eta^3 \mp \tfrac2r \cos\phi\, \eta^4.$$
Once we restrict to $M$, then $\eta^4=\kappa^\pm_3=0$, so that
$$
d\bg^+ = -\frac1r \cos \phi\, \eta^3 \bg^+
+\frac1r  (\eta^0 - \sin \phi\, \eta^3)\ri \bg^+ \quad \mod \kappa^+_1, \kappa^+_2.
$$
Thus, along a curve in $M$ tangent to the leaves of $\chi^+$, the derivative of $\bg^+$
is a complex multiple of the value of $\bg^+$.  Therefore, the complex line spanned
by $\bg^+$ remains fixed as we move along the leaves of $\chi^+$.
\end{proof}

\subsection{Proof of Main Result}\label{mainproof}
Given any point of $M$, we can construct an adapted lift $f$ on a neighborhood
of that point.  To prove the first part of Assertion 1 in the theorem, we need to to show that
the composition $\bg^\pm_\C\circ f$ is unchanged when we modify the adapted lift.
For simplicity, we will show that $\bg^+_\C \circ f$ is unchanged, the
argument for $\bg^-_\C$ being similar.

At a given point of $M$, any two adapted lifts differ either by
moving the point $\bz$ along the fiber of the quotient map $\pi:\HH \to \CH^2$,
or by rotating $\eg_1$ and $\eg_2$ while keeping their span fixed.
We will define vector fields on $G$ whose trajectories correspond to
modifying the frame in these ways.  Note that by Lemma \ref{IIlemma}, an adapted lift of a hypersurface
is a 3-dimensional submanifold $\Sigma \subset G$ on which the 1-form $\eta^4$ pulls back
to be zero.  (As well, because the projection of this submanifold to $\HH$ must be transverse
to the fibers of $\pi$, the 3-form $\eta^1 \& \eta^2 \& \eta^3$ must pull back to be
nonzero at each point of $\Sigma$.)

First,
a left-invariant vector field $X$ on $G$ that moves $\bz$ in the direction of $\eg_0$
(i.e., tangent to the fiber of $\pi$) must, by \eqref{dzed}, have the property that $X \intprod \eta^0$
is a nonzero constant, while its interior product with $\eta^1, \ldots, \eta^4$ must be zero.
In order that motion along the trajectories of $X$ take adapted lifts to adapted lifts,
we need $\Lie_X \eta^4$ to be a multiple of $\eta^4$.  Because
$$d\eta^4 = -\w^4_1 \& \eta^1 - \w^4_2 \& \eta^2 - \w^4_3 \& \eta^3,$$
the interior product of $X$ with $\w^4_1, \w^4_2, \w^4_3$ must be zero.
(To define $X$ uniquely, we normalize $X \intprod \eta^0=1$ and require also
that $X \intprod \w^2_1=0$.)
Then from \eqref{dgplus} we see that
$$X \intprod d\bg^+ = \frac{\ri}r \bg^+.$$
Thus, this change of adapted lift only changes $\bg^+$ by multiplication by
a unit modulus constant.

Next, note that \eqref{dhatalpha} specializes to
$$d\eg_1 = \w^2_1 \eg_2 + \w^3_1 \eg_3 + \w^4_1 \eg_4 + \tfrac1{r^2}\eta^1 \zz -\tfrac1r \eta^2 \eg_0.$$
Therefore, a left-invariant vector field $Y$ on $G$ that corresponds to rotating
$\eg_1$ in the direction of $\eg_2$ satisfies
$Y \intprod \w^2_1=0$, while the interior product of $Y$ with all the
other 1-forms $\eta^0, \ldots, \eta^4$, $\eta^4_1$, $\eta^4_2$ and $\eta^4_3$ is zero.
(To see that motion along the trajectories of $Y$ takes adapted lifts to adapted lifts,
it is easy to check that $\Lie_Y \eta^4=0$.)  Then we easily calculate that
$Y \intprod d\bg^+=0$.

It follows from Proposition \ref{intprop} that the images of $\sigma^\pm = \bg^\pm_\C \circ f$
are curves.  (Moreover, because for each choice of sign, one of $\kappa^\pm_1, \kappa^\pm_2$ is nonzero
at every point of $M$, \eqref{dgplus} shows that the images are regular curves.)
To prove the rest of Assertion 1, we need to show that these are contact curves.

Recall that the maps $\bg^\pm_\C$ take value in the projectivization of $\pi(\V)$ of the null
cone $\V \subset \C^3$.  As explained in \S\ref{gausschar}, this may be identified with $S^3$.
We now give a characterization of the standard contact structure in $S^3$ under this identification:

\begin{lemma}\label{contactcond}  Let $\mu: I \to \pi(\V)$ be a regular curve (where $I$ is an interval on
the real line), and let $\bn:I \to \V$ be any lift of $\mu$.  Then the image of $\mu$
is a contact curve if and only if
\begin{equation}\label{nullcond}
\left\langle \dfrac{d\bn}{d t},\ri \bn(t)\right\rangle =0.
\end{equation}
\end{lemma}
\begin{proof}
Note that condition \eqref{nullcond} is independent of choice of lift $\bn$.

Suppose that $\bn(t)=(z_0(t),z_1(t),z_2(t))$ with $|z_0|^2=|z_1|^2 + |z_2|^2$.
Then
$$\left\langle \dfrac{d\bn}{dt} ,  \ri  \bn \right\rangle
=\realpart\left(-\ri\left(
\overline{z_1} z_1' +\overline{z_2}z_2' -\overline{z_0} z_0'\right)\right).$$
Meanwhile, the corresponding curve in $\C^2$, given by
$w_1=z_1/z_0$ and $w_2 = z_2/z_0$, satisfies
\begin{align*}\left\langle \dfrac{d\bw}{dt}, \ri \bw \right\rangle
&= \realpart\left( \dfrac{-\ri}{z_0 |z_0|^2}
\left( \overline{z_1}(z_0 z_1'-z_1 z_0')
+ \overline{z_2}(z_0 z_2'-z_2 z_0')\right) \right)\\
&= \realpart\left( \dfrac{-\ri}{|z_0|^2}
\left(\overline{z_1} z_1' +\overline{z_2}z_2' -\overline{z_0} z_0'\right)
\right),
\end{align*}
where we use the usual Hermitian inner product on $\C^2$.
(Note that at $w\in S^3$, the vector $\ri w$ is tangent to the fibers of
the Hopf fibration.)
Therefore, $\mu(t)$ is a contact curve in the space of complex null lines if and only
if $\bw(t)$ is a contact curve in $S^3 \subset \C^2$.
\end{proof}

We finish the proof of Assertion 1 with the following lemma, which will
also be used in the proof of Assertion 2.

\begin{lemma}\label{mukappa} Let $\mu:I \to S^3$ be any regular curve (where $I$ is an interval
on the real line), and let $\gamma:I \to G$ be
any lift such that $\bg^\pm_\C \circ \gamma = \mu$.  Then $\mu$ is a contact curve if and only
if $\gamma^* \kappa^\pm_3 =0$.
\end{lemma}
\begin{proof}We will verify the assertion for $\bg^+_\C$, the argument for $\bg^-_\C$ being similar.
First, note that $\bg^+ \circ \gamma$ is a lift of $\mu$ into the null cone $\V$.
Let $t$ be a coordinate on $I$, and let $\bg^+(t)$ stand for $\bg^+(\gamma(t))$.
Then, applying Lemma \ref{contactcond} and using \eqref{dgplus}, we
see that $\mu$ is a contact curve if and only if
$$0 = \left\langle \dfrac{d\bg^+}{d t}, \ri \bg^+(t)\right\rangle
= \dfrac{\di}{\di t} \intprod (\gamma^* \kappa^+_3).$$
\end{proof}

To begin the proof of Assertion 2, note that the nine 1-forms $\eta^0, \eta^3, \w^1_2$
and $\kappa^\pm_1$, $\kappa^\pm_2$, $\kappa^\pm_3$ comprise a coframe on $G$.
Thus, by \eqref{dgplus}, the map $\bg^+_\C$ has constant rank 3 (and the same is true for $\bg^-_\C$).
Suppose that $\CC_1$ is an embedded contact curve in $S^3$, and let $N\subset G$ be its inverse
image under $\bg^+_\C$.  Then $N$ is a 7-dimensional submanifold.  Because
$\kappa^+_3$ must pull back to be zero on $N$, and $\kappa^+_1, \kappa^+_2$ pull back to be linearly
dependent on $N$, then the remaining forms $\eta^0, \eta^3, \w^1_2, \kappa^-_1, \kappa^-_2, \kappa^-_3$
pull back to $N$ to be linearly independent.  In particular, the restriction of $\bg^-_\C$ to $N$ still has
rank 3.  It follows that the intersection $P$ of the inverse images of the contact curves
is a smooth submanifold of dimension 5.

Since the maps $\bg^\pm_\C$ are constant along the trajectories
of $X$ and $Y$, then $P$ is foliated by these trajectories.  It follows that the image $M=\rho(P)$ inside $\CH^2$ is three-dimensional.  Moreover, any adapted lift of $M$ lies inside $P$.  Because
the 1-forms $\kappa^+_3$ and $\kappa^-_3$ pull back to be zero on $P$, the same is true
for $\w^4_3 - \alpha \eta^3$.  So, for any adapted lift $f$ of $M$,
$f^*(\w^4_3 -\alpha \eta^3)=0$.  It then follows from Lemma \ref{IIlemma} that $M$ is a Hopf hypersurface
with Hopf principal curvature $\alpha$.

This finishes the proof of Theorem \ref{mainresult}.

\subsection{Borderline Case}\label{sideproof}

In the borderline case, where $\alpha =\pm 2/r$
the functions $\bg^+$ and $\bg^-$ coincide.  However, we can still use this modified Gauss map to associate
a single contact curve to a Hopf hypersurface, and to some extent we can reconstruct the hypersurface
from the curve.  For simplicity, say that $\alpha=2/r$, whence $\sin\phi=1$ and $\cos\phi=0$.
Let
$$\bg = \eg_0 - \eg_3,$$
and let $\bg_\C = \pi \circ \bg$ denote the projectivization of $\bg$.
Then \eqref{dgplus} specializes to
$$d\bg =\frac1r \eta^4 \bg
+\frac1r  (\eta^0 + \eta^3)\ri\bg
-\kappa_1 \eg_2 + \kappa_2 \eg_1 -\kappa_3 \eg_4,
$$
where we define
$$
\kappa_1 = \w^4_1 -\frac1r \eta^1, \quad
\kappa_2 = \w^4_2 -\frac1r \eta^2,\quad
\kappa_3 = \w^4_3 - \tfrac{2}{r} \eta^3= \w^4_3 - \alpha \eta^3.
$$
The proof in \S\ref{mainproof} for Assertion 1 may be suitably modified to show that
the composition $\sigma = \bg_\C \circ f$, where $f$ is an arbitrary adapted lift of a Hopf hypersurface
$M$ with Hopf principal curvature $\alpha$, gives a well-defined map to $S^3$, and
whose image is a contact curve.
(In fact, the proof of Lemma \ref{mukappa} can
be modified to show that a curve $\mu:I\to S^3$ is contact if and only if $\gamma^* \kappa_3 = 0$,
where $\gamma: I \to G$ is any lift such that $\bg_\C\circ \gamma =\mu$.)
Note that \eqref{splitsig} implies that, at each point of $M$, one of the other principal curvatures
$\lambda,\nu$ must equal $1/r$.  If they are not both equal to $1/r$, then one of $\kappa_1, \kappa_2$
is nonzero at each point, and rank of $\sigma$ is one.

Thus, given an arbitrary regular contact curve $\CC$ in $S^3$, the inverse image $N$ under $\bg_\C$ of this
curve is a codimension two submanifold of $G$ on which $\kappa_3$ pulls back to be zero
and $\kappa_1, \kappa_2$ pull back to be linearly dependent.  We wish to construct a three-dimensional
submanifold $\Sigma \subset N$ which will be an adapted lift of a Hopf hypersurface.  To do this,
we only need to ensure that $\eta^4$ pulls back to be zero on $\Sigma$,
and the forms $\eta^1, \eta^2, \eta^3$ pull back to be linearly independent.  (We refer to
the latter as the {\em independence condition} for $\Sigma$.)
Then the vanishing of $\kappa_3$ will ensure that the Hopf condition holds.

It is easy to check that
$$d\eta^4 \equiv -\kappa_1 \& \eta^1 -\kappa_2 \& \eta^2 \quad \mod \kappa_3.$$
Thus, because the span of \{$\kappa_1, \kappa_2$\} is one-dimensional on $N$ (the regularity
of the contact curve rules out the vanishing of both $\kappa_1$ and $\kappa_2$ at any point),
the Pfaff rank of $\eta^4$ is one.  (In other words, $\eta^4 \& d \eta^4 \ne 0$
but $\eta^4 \& d \eta^4 \& d\eta^4 = 0$.)  It follows by Pfaff's theorem
(see \cite{cfb} Chapter One) that around any point of $N$ there is an open set in $N$ and
a local coordinate
system $x_0, \ldots, x_6$ defined on that set, in which $\eta^4$ is a multiple of $dx_1 - x_2 dx_0$.
It follows that by setting $x_1 = \psi(x_0)$, $x_2 =\psi'(x_0)$ for an arbitrary function $\psi$, we
can construct a five-dimensional submanifold $P$ on which both $\kappa_3$ and $\eta^4$ pull back to be zero.
As argued in \S\ref{mainproof}, the trajectories of $X$ and $Y$ foliate $P$,
and its image under $\rho$ is a Hopf hypersurface in $\CH^2$.

To summarize, here is the analogue of Theorem \ref{mainresult} in the borderline case:
\begin{theorem}\label{sideresult}\
\begin{enumerate}
\item Let $M$ be an oriented Hopf hypersurface in $\CH^2$
with Hopf principal curvature $\alpha=\pm2/r$.  Then there is a well-defined
map $\sigma:M \to S^3$, such that $\sigma = \bg_\C \circ f$ for any adapted
lift $f$ of $M$, and such that the image $\sigma(M)$ is a contact curve in $S^3$.
(The map $\sigma$ is regular at  points where the principal curvatures of $M$ are distinct.)

\item Conversely, given any regular contact curve $\CC$ in $S^3$, let $N\subset G$ be
its inverse image under $\bg_\C$.  Given $p\in N$, there exists an open neighborhood
$U \subset N$ containing $p$, and a family of five-dimensional manifolds $P \subset U$
containing $p$, such that $\pi(P)$ is a Hopf hypersurface with principal curvature
$\pm 2/r$.  (The family is parametrized by a choice of one function of one real variable.)
\end{enumerate}
\end{theorem}

\section{Afterword}
Our original approach to the construction of Hopf hypersurfaces in $\CH^2$ was to analyze
the Hopf condition using the techniques of exterior differential systems (EDS).
From this point of view, adapted lifts of Hopf hypersurfaces are 3-dimensional integral
submanifolds of the EDS on $G$ generated by the 1-forms $\eta^4$ and
$\w^4_3 - \alpha \eta^3$.
Applying Cartan's Test shows that this Pfaffian system is involutive with last nonzero
character $s_1=2$.  This indicates that locally-defined solutions may be constructed by solving
a sequence of Cauchy problems,
in which the last time we have freedom to specify the initial data,
that choice amounts to
two functions of one variable.
The existence of these solutions is guaranteed by
the Cartan-K\"ahler theorem.

However, the Cartan-K\"ahler theorem does not give an explicit construction
for solutions, so one has to work harder.
In fact, our analysis shows that, when $\alpha$ is strictly between $-2/r$ and $2/r$,
this EDS is hyperbolic, and is integrable by the method of Darboux.  The construction
scheme for Hopf hypersurfaces given by Theorem \ref{mainresult} is based on this observation.
We note that the function count given by Cartan's Test is explicitly realized by this construction,
since specifying a contact curve in a three-dimensional manifold amounts to choosing
an arbitrary function of one variable (because the contact form may be written as $dy - z\,dx$
in suitable local coordinates).

\end{document}